\documentclass{amsart}

\usepackage{amsmath,amsthm,amsfonts,amscd,amssymb}
\usepackage[pdfborder={0 0 0}]{hyperref}
\usepackage{enumerate}

\newcommand{\ZZ}{\mathbb{Z}}

\newtheorem{thm}{Theorem}[section]

\newtheorem{lem}[thm]{Lemma}

\theoremstyle{definition}

\theoremstyle{remark}
\newtheorem{rem}{Remark}[section]

\begin{document}

\title{Well-rounded zeta-function of planar arithmetic lattices}
\author{Lenny Fukshansky}\thanks{The author was partially supported by a grant from the Simons Foundation (\#208969 to Lenny Fukshansky) and by the NSA Young Investigator Grant \#1210223.}

\address{Department of Mathematics, 850 Columbia Avenue, Claremont McKenna College, Claremont, CA 91711}
\email{lenny@cmc.edu}



\begin{abstract}
We investigate the properties of the zeta-function of well-rounded sublattices of a fixed arithmetic lattice in the plane. In particular, we show that this function has abscissa of convergence at $s=1$ with a real pole of order 2, improving upon a result of \cite{kuehnlein}. We use this result to show that the number of well-rounded sublattices of a planar arithmetic lattice of index less or equal $N$ is $O(N \log N)$ as $N \to \infty$. To obtain these results, we produce a description of integral well-rounded sublattices of a fixed planar integral well-rounded lattice and investigate convergence properties of a zeta-function of similarity classes of such lattices, building on the results of \cite{fletcher_jones}.
\end{abstract}

\maketitle

\def\A{{\mathcal A}}
\def\AA{{\mathfrak A}}
\def\B{{\mathcal B}}
\def\C{{\mathcal C}}
\def\D{{\mathcal D}}
\def\EE{{\mathfrak E}}
\def\F{{\mathcal F}}
\def\x{{\mathcal H}}
\def\I{{\mathcal I}}
\def\II{{\mathfrak I}}
\def\J{{\mathcal J}}
\def\K{{\mathcal K}}
\def\kk{{\mathfrak K}}
\def\L{{\mathcal L}}
\def\LL{{\mathfrak L}}
\def\M{{\mathcal M}}
\def\mm{{\mathfrak m}}
\def\MM{{\mathfrak M}}
\def\N{{\mathcal N}}
\def\O{{\mathcal O}}
\def\OO{{\mathfrak O}}
\def\PP{{\mathfrak P}}
\def\R{{\mathcal R}}
\def\PNR{{\mathcal P_N(\real)}}
\def\PMNR{{\mathcal P^M_N(\real)}}
\def\PdNR{{\mathcal P^d_N(\real)}}
\def\s{{\mathcal S}}
\def\V{{\mathcal V}}
\def\X{{\mathcal X}}
\def\Y{{\mathcal Y}}
\def\Z{{\mathcal Z}}
\def\H{{\mathcal H}}
\def\cee{{\mathbb C}}
\def\Nn{{\mathbb N}}
\def\pee{{\mathbb P}}
\def\que{{\mathbb Q}}
\def\QQ{{\mathbb Q}}
\def\real{{\mathbb R}}
\def\RR{{\mathbb R}}
\def\zed{{\mathbb Z}}
\def\ZZ{{\mathbb Z}}
\def\aaa{{\mathbb A}}
\def\ff{{\mathbb F}}
\def\HDelta{{\it \Delta}}
\def\kk{{\mathfrak K}}
\def\qbar{{\overline{\mathbb Q}}}
\def\kbar{{\overline{K}}}
\def\ybar{{\overline{Y}}}
\def\kkbar{{\overline{\mathfrak K}}}
\def\ubar{{\overline{U}}}
\def\eps{{\varepsilon}}
\def\ahat{{\hat \alpha}}
\def\bhat{{\hat \beta}}
\def\gt{{\tilde \gamma}}
\def\h{{\tfrac12}}
\def\be{{\boldsymbol e}}
\def\bei{{\boldsymbol e_i}}
\def\bc{{\boldsymbol c}}
\def\bm{{\boldsymbol m}}
\def\bk{{\boldsymbol k}}
\def\bi{{\boldsymbol i}}
\def\bl{{\boldsymbol l}}
\def\bq{{\boldsymbol q}}
\def\bu{{\boldsymbol u}}
\def\bt{{\boldsymbol t}}
\def\bs{{\boldsymbol s}}
\def\bv{{\boldsymbol v}}
\def\bw{{\boldsymbol w}}
\def\bx{{\boldsymbol x}}
\def\bX{{\boldsymbol X}}
\def\bz{{\boldsymbol z}}
\def\bwy{{\boldsymbol y}}
\def\bY{{\boldsymbol Y}}
\def\bL{{\boldsymbol L}}
\def\ba{{\boldsymbol a}}
\def\bb{{\boldsymbol b}}
\def\bet{{\boldsymbol\eta}}
\def\bxi{{\boldsymbol\xi}}
\def\bo{{\boldsymbol 0}}
\def\bone{{\boldsymbol 1}}
\def\bol{{\boldsymbol 1}_L}
\def\ep{\varepsilon}
\def\p{\boldsymbol\varphi}
\def\q{\boldsymbol\psi}
\def\rank{\operatorname{rank}}
\def\aut{\operatorname{Aut}}
\def\lcm{\operatorname{lcm}}
\def\sgn{\operatorname{sgn}}
\def\spn{\operatorname{span}}
\def\md{\operatorname{mod}}
\def\Norm{\operatorname{Norm}}
\def\dim{\operatorname{dim}}
\def\det{\operatorname{det}}
\def\Vol{\operatorname{Vol}}
\def\rk{\operatorname{rk}}
\def\ord{\operatorname{ord}}
\def\ker{\operatorname{ker}}
\def\div{\operatorname{div}}
\def\Gal{\operatorname{Gal}}
\def\GL{\operatorname{GL}}
\def\SL{\operatorname{SL}}
\def\SNR{\operatorname{SNR}}
\def\WR{\operatorname{WR}}
\def\IWR{\operatorname{IWR}}
\def\scg{\operatorname{\left< \Gamma \right>}}
\def\swrh{\operatorname{Sim_{WR}(\Lambda_h)}}
\def\ch{\operatorname{C_h}}
\def\cht{\operatorname{C_h(\theta)}}
\def\scgt{\operatorname{\left< \Gamma_{\theta} \right>}}
\def\scgmn{\operatorname{\left< \Gamma_{m,n} \right>}}
\def\gat{\operatorname{\Omega_{\theta}}}
\def\mn{\operatorname{mn}}
\def\disc{\operatorname{disc}}
\def\Re{\operatorname{Re}}
\def\lcm{\operatorname{lcm}}

\section{Introduction}
\label{intro}

Let $\Lambda = A\zed^2  \subset \real^2$ be a lattice of full rank in the plane, where $A=(\ba_1 \ba_2)$ is a basis matrix. The corresponding norm form is defined as 
$$Q_A(\bx) = \bx^t A^t A \bx.$$
We say that $\Lambda$ is {\it arithmetic} if the entries of the matrix $A^tA$ generate a 1-dimensional $\que$-vector subspace of $\real$. This property is easily seen to be independent of the choice of a basis. We define $\det(\Lambda)$ to be $|\det(A)|$, again independent of the basis choice, and (squared) {\it minimum} or {\it minimal norm}
$$|\Lambda| = \min \{ \|\bx\|^2 : \bx \in \Lambda \setminus \{\bo\} \} = \min \{ Q_A(\bwy) : \bwy \in \zed^2 \setminus \{\bo\} \},$$
where $\|\ \|$ stands for the usual Euclidean norm. Then each $\bx \in \Lambda$ such that $\|\bx\|^2 = |\Lambda|$ is called a {\it minimal vector}, and the set of minimal vectors of $\Lambda$ is denoted by $S(\Lambda)$. A planar lattice $\Lambda$ is called {\it well-rounded} (abbreviated WR) if the set $S(\Lambda)$ contains a basis for $\Lambda$; we will refer to such a basis as a {\it minimal basis} for~$\Lambda$.

While in this note we focus on the planar case, the notion of WR lattices is defined in every dimension: a full-rank lattice in $\real^N$ is WR if it contains $N$ linearly independent minimal vectors -- the fact that these form a basis for the lattice is a low-dimensional phenomenon, only valid for $N \leq 4$. WR lattices are important in discrete optimization, in particular in the investigation of sphere packing, sphere covering, and kissing number problems \cite{martinet}, as well as in coding theory \cite{esm}. Properties of WR lattices have also been investigated in \cite{mcmullen} in connection with Minkowski's conjecture and in \cite{lf:robins} in connection with the linear Diophantine problem of Frobenius. Furthermore, WR lattices are used in cohomology computations of $\SL_N(\zed)$ and its subgroups \cite{ash}. These considerations motivate the study of distribution properties of WR lattices. Distribution of WR lattices in the plane has been studied in \cite{wr1}, \cite{wr2}, \cite{wr3}, \cite{fletcher_jones}, \cite{kuehnlein}. In particular, these papers investigate various aspects of distribution properties of WR sublattices of a fixed planar lattice.

An important equivalence relation on lattices is geometric similarity: two lattices $\Lambda_1, \Lambda_2 \subset \real^2$ are called {\it similar}, denoted $\Lambda_1 \sim \Lambda_2$, if there exists  $\alpha \in \real_{>0}$ and $U \in O_2(\real)$ such that $\Lambda_2 = \alpha U \Lambda_1$. It is easy to see that similar lattices have the same algebraic  structure, i.e., for every sublattice $\Gamma_1$ of a fixed index in $\Lambda_1$ there is a sublattice $\Gamma_2$ of the same index in $\Lambda_2$ so that $\Gamma_1 \sim \Gamma_2$. A WR lattice can only be similar to another WR lattice, so it makes sense to speak of WR similarity classes of lattices. In \cite{kuehnlein} it has been proved that a planar lattice contains infinitely many non-similar WR sublattices if and only if it contains one. This is always the case for arithmetic planar lattices. If the lattice in question is not arithmetic, it may still have infinitely many non-similar WR sublattices depending on the value of a certain invariant described in \cite{kuehnlein}. In any case, it appears that non-arithmetic planar lattices contain fewer WR sublattices than arithmetic ones in the sense which we discuss below.

Given an infinite finitely generated group $G$, it is a much-studied problem to determine the asymptotic growth of $\# \left\{ H \leq G : \left| G:H \right| \leq N \right\}$, the number of subgroups of index no greater than $N$, as $N \to \infty$ (see \cite{lubot}). One approach that has been used by different authors with great success entails looking at the analytic properties of the corresponding Dirichlet-series generating function $\sum_{H \leq G} \left| G:H \right|^{-s}$ and then using some Tauberian theorem to deduce information about the rate of growth of partial sums of its coefficients (see \cite{sautoy}, as well as Chapter~15 of \cite{lubot}). In case $G$ is a free abelian group of rank 2, i.e. a planar lattice, this Dirichlet series allows to count sublattices of finite index, and is a particular instance of the Solomon zeta-function (see \cite{reiner}, \cite{solomon}). We will use a similar approach while restricting to just WR sublattices, which is a more delicate arithmetic problem.

Fix a planar lattice $\Omega$, and define the {\it zeta-function of WR sublattices} of $\Omega$ to be
$$\zeta_{\WR}(\Omega,s) = \sum_{\WR \Lambda \subseteq \Omega} \frac{1}{\left| \Omega : \Lambda \right|^s} = \sum_{n=1}^{\infty} \frac{\# \{\WR \Lambda \subseteq \Omega : \left| \Omega : \Lambda \right| = n\}}{n^s}$$
for $s \in \cee$. The rate of growth of coefficients of this function can be conveyed by studying its abscissa of convergence and behavior of the function near it. For brevity of notation, we will say that an arbitrary Dirichlet series $f(s) = \sum_{n=1}^{\infty} a_n n^{-s}$ has an {\it abscissa of convergence} with a {\it real pole of order $\mu$} at $s=\rho$ if $f(s)$ is absolutely convergent for $\Re(s) > \rho$, and for $s \in \real$
\begin{equation}
\label{lim_def}
\lim_{s \to \rho^+} (s-\rho)^{\mu} \sum_{n=1}^{\infty} \frac{a_n}{n^s}
\end{equation}
exists and is nonzero. Notice that this notion does not imply existence of analytic continuation for $f(s)$, but is merely a statement about the rate of growth of the coefficients of $f(s)$, which is precisely what we require. For instance, in \cite{wr1} and \cite{wr2} it has been established that $\zeta_{\WR}(\zed^2,s)$ has abscissa of convergence with a real pole of order 2 at $s=1$. Furthermore, it has been shown in \cite{kuehnlein} that if $\Omega$ is a non-arithmetic planar lattice containing WR sublattices, then $\zeta_{\WR}(\Omega,s)$ has abscissa of convergence with a real pole of order 1 at $s=1$ (in fact, Lemma~3.3 of~\cite{kuehnlein} combined with Theorem~4 on p.~158 of~\cite{lang} imply the existence of analytic continuation of $\zeta_{\WR}(\Omega,s)$ in this situation to $\Re(s) > 1-\eps$ for some $\eps > 0$ with a pole of order 1 at $s=1$). It is natural to expect that the situation for any arithmetic lattice is the same as it is for $\zed^2$; in fact, another result of \cite{kuehnlein} states that for any arithmetic lattice $\Omega$, $\zeta_{\WR}(\Omega,s)$ has abscissa of convergence at $s=1$, and it is conjectured that it has a pole of order 2 at $s=1$. The main goal of the present paper is to prove the following result in this direction.

\begin{thm} \label{main} Let $\Omega$ be a planar arithmetic lattice. Then $\zeta_{\WR}(\Omega,s)$ has abscissa of convergence with a real pole of order 2 at $s=1$ in the sense of \eqref{lim_def} above. Moreover,
\begin{equation}
\label{growth_bnd_1}
\# \{\WR\ \Lambda \subseteq \Omega : \left| \Omega : \Lambda \right| \leq N\} = O(N \log N)
\end{equation}
as $N \to \infty$.
\end{thm}
\smallskip

\begin{rem} \label{WR_growth} To compare, Theorem~4.20 of \cite{sautoy} combined with Lemma~3.3 (and the Corollary following it) of \cite{kuehnlein} imply that if $\Omega$ is a non-arithmetic planar lattice containing WR sublattices, then the right hand side of \eqref{growth_bnd_1} is equal to $O(N)$. It should be pointed out that by writing that a function of $N$ is equal to $O(N \log N)$ (respectively, $O(N)$) we mean here that it is asymptotically bounded from above and below by nonzero multiples of $N \log N$ (respectively,~$N$). On the other hand, it is a well known fact (outlined, for example, on p. 793 of \cite{sautoy}) that for any planar lattice~$\Omega$,
\begin{equation}
\label{growth_bnd_2}
\# \{ \Lambda \subseteq \Omega : \left| \Omega : \Lambda \right| \leq N\} \sim \left( \pi^2/12 \right) N^2
\end{equation}
as $N \to \infty$.
\end{rem}

The organization of this paper is as follows. In Section~\ref{IWR}, we start by reducing the problem to integral WR (abbreviated IWR) lattices in Lemma~\ref{reduce}: a planar lattice $\Lambda = A\zed^2$ is called {\it integral} if the coefficient matrix $A^t A$ of its quadratic form $Q_A$ has integer entries (this definition does not depend on the choice of a basis).  We then introduce zeta-functions of similarity classes of planar IWR lattices, objects of independent interest, and study their convergence properties in Theorem~\ref{IWR_zeta}. Our arguments build on the parameterization of planar IWR lattices obtained in \cite{fletcher_jones}. In Section~\ref{IWR_subl} we continue using this parameterization to obtain an explicit description of IWR sublattices of a fixed planar IWR lattice which are similar to another fixed IWR lattice (Theorem~\ref{two_IWR_1}), and use it to determine convergence properties of the Dirichlet series generating function of all such sublattices (Lemma~\ref{zeta_two}). Finally, in Lemma~\ref{zeta_two_1} we decompose $\zeta_{\WR}(\Omega,s)$ for a fixed IWR planar lattice $\Omega$ into a sum over similarity classes of sublattices and observe that this sum can be represented as a product of the two different types of Dirichlet series that we investigated above; hence the result of Theorem~\ref{main} follows by Lemma~\ref{reduce}.
\bigskip

\section{Integral WR lattices in the plane}
\label{IWR}

Integral lattices are central objects in arithmetic theory of quadratic forms and in lattice theory. IWR lattices  have recently been studied in \cite{fletcher_jones}. The significance of IWR planar lattices for our purposes is reflected in the following reduction lemma.

\begin{lem} \label{reduce} Let $\Omega$ be an arithmetic planar lattice. Then there exists some IWR planar lattice $\Lambda$ such that $\zeta_{\WR}(\Omega,s)$ has the same abscissa of convergence with pole of the same order as $\zeta_{\WR}(\Lambda,s)$.
\end{lem}

\proof Lemma 2.1 of \cite{kuehnlein} guarantees that $\Omega$ has a WR sublattice, call it $\Omega'$; naturally, $\Omega'$ must also be arithmetic. Let $A$ be a basis matrix for $\Omega'$, then entries of $A^tA$ span a 1-dimensional vector space over $\que$, meaning that there exists $\alpha \in \real_{>0}$ such that the matrix $\alpha A^tA$ is integral. Then the lattice $\Lambda := \sqrt{\alpha} A \zed^2$ is integral and is similar to $\Omega'$, hence it is also WR. Since $\Lambda$ is just a scalar multiple of $\Omega'$, it is clear that $\zeta_{\WR}(\Lambda,s)$ has the same abscissa of convergence with pole of the same order as $\zeta_{\WR}(\Omega',s)$, which is the same as that of $\zeta_{\WR}(\Omega,s)$ by Lemma~3.2 of~\cite{kuehnlein}.
\endproof

Moreover, it is easy to see that these properties of zeta-function of WR sublattices are preserved under similarity.

\begin{lem} \label{sim_red} Assume that $\Lambda_1, \Lambda_2$ are two planar lattices such that $\Lambda_1 \sim \Lambda_2$. Then $\zeta_{\WR}(\Lambda_1,s) = \zeta_{\WR}(\Lambda_2,s)$.
\end{lem}

\proof Similar lattices have the same numbers of WR sublattices of the same indices. The statement of the lemma follows immediately. 
\endproof

Lemmas~\ref{reduce} and~\ref{sim_red} imply that we can focus our attention on similarity classes of IWR lattices to prove Theorem~\ref{main}. Integrality is not preserved under similarity, however a WR similarity class may or may not contain integral lattices. WR similarity classes containing integral lattices, we will call them IWR similarity classes, have been studied in \cite{fletcher_jones} -- these are precisely the WR similarity classes containing arithmetic lattices. Let us write $\left< \Lambda \right>$ for the similarity class of the lattice $\Lambda$, then a result of \cite{fletcher_jones} states that the set of IWR similarity classes is
$$\IWR = \left\{ \left< \Gamma_D(p,q) \right> : \Gamma_D(p,q) = \frac{1}{\sqrt{q}} \begin{pmatrix} q & p \\ 0 & r\sqrt{D} \end{pmatrix} \zed^2 \right\},$$
where $(p,r,q,D)$ are all positive integer 4-tuples satisfying
\begin{equation}
\label{prqD}
p^2+Dr^2=q^2,\ \gcd(p,q)=1,\ \frac{p}{q} \leq \frac{1}{2}, \text{ and } D \text{ squarefree}.
\end{equation}
It is also discussed in \cite{fletcher_jones} that $\Gamma_D(p,q)$ is a {\it minimal} integral lattice with respect to norm in its similarity class. In particular, every integral lattice $\Lambda \in \left< \Gamma_D(p,q) \right>$ is of the form $\Lambda = \sqrt{k}\ U \Gamma_D(p,q)$ for some $k \in \zed_{>0}$, $U \in O_2(\real)$, and so
$$|\Lambda| \geq |\Gamma_D(p,q)| = q.$$
The set $\IWR$ can be represented as
$$\IWR = \bigsqcup_{D \in \zed_{>0} \text{ squarefree}} \IWR(D),$$
where for each fixed positive squarefree integer $D$, $\IWR(D) := \left\{ \left< \Gamma_D(p,q) \right> \right\}$ is the set of IWR similarity classes of {\it type} $D$.

Let us define the {\it minimum} and {\it determinant zeta-functions} of IWR similarity classes of type $D$ in the plane:
\begin{equation}
\label{z^m_IWR}
\zeta^m_{\IWR(D)}(s) = \sum_{\left< \Gamma_D(p,q) \right> \in \IWR(D)} \frac{1}{|\Gamma_D(p,q)|^s} = \sum_{\left< \Gamma_D(p,q) \right> \in \IWR(D)} \frac{1}{q^s},
\end{equation}
and
\begin{equation}
\label{z^d_IWR}
\zeta^d_{\IWR(D)}(s) = \sum_{\left< \Gamma_D(p,q) \right> \in \IWR(D)} \frac{1}{\det \Gamma_D(p,q)^s} = \frac{1}{D^{s/2}} \sum_{\left< \Gamma_D(p,q) \right> \in \IWR(D)} \frac{1}{r^s},
\end{equation}
where $s \in \cee$. Since
\begin{equation}
\label{rq}
\frac{\sqrt{3}}{2} \times \frac{1}{\sqrt{D}} \times q \leq r \leq \frac{1}{\sqrt{D}} \times q,
\end{equation}
we have
\begin{equation}
\label{zeta_md_ineq}
\zeta^m_{\IWR(D)}(s) \leq \zeta^d_{\IWR(D)}(s) \leq \left( \frac{2}{\sqrt{3}} \right)^s \zeta^m_{\IWR(D)}(s)
\end{equation}
for all real $s$, and so $\zeta^m_{\IWR(D)}(s)$ and $\zeta^d_{\IWR(D)}(s)$ have the same convergence properties. We can establish the following result.

\begin{thm} \label{IWR_zeta} For every real value of $s > 1$,
\begin{equation}
\label{zeta_bnd}
\frac{1}{\left( 2\sqrt{3D} \right)^s} \frac{\zeta(2s-1)}{\zeta(2s)} \leq \zeta^d_{\IWR(D)}(s)  \leq \left( \frac{2}{\sqrt{3}} \right)^s \zeta^m_{\IWR(D)}(s) \leq \left( \frac{4D}{\sqrt{3}} \right)^s \zeta_{\que(\sqrt{-D})}(s),
\end{equation}
where $\zeta(s)$ is the Riemann zeta-function and $\zeta_{\que(\sqrt{-D})}(s)$ is the Dedekind zeta-function of the imaginary quadratic number field $\que(\sqrt{-D})$. Hence the Dirichlet series $\zeta^d_{\IWR(D)}(s)$ and $\zeta^m_{\IWR(D)}(s)$ are absolutely convergent for $\Re(s) > 1$, and for $s \in \real$ the limits
\begin{equation}
\label{lim_md}
\lim_{s \to 1^+} (s-1) \zeta^d_{\IWR(D)}(s),\ \lim_{s \to 1^+} (s-1) \zeta^m_{\IWR(D)}(s)
\end{equation}
exist and are nonzero. Moreover, the $N$-th partial sums of coefficients of these Dirichlet series are equal to $O(N)$ as $N \to \infty$.
\end{thm}

\proof Let $D$ be a fixed positive squarefree integer. Lemma~1.3 of \cite{fletcher_jones} guarantees that  $p,r,q \in \zed_{>0}$ satisfy \eqref{prqD} if and only if
\begin{equation}
\label{mn_par}
p =  \frac{| m^2-Dn^2 |}{2^e \gcd(m,D)},\ r = \frac{2mn}{2^e \gcd(m,D)},\ q = \frac{m^2 + Dn^2}{2^e \gcd(m,D)},
\end{equation}
for some $m,n \in \zed$ with $\gcd(m,n)=1$ and $\sqrt{\frac{D}{3}} \leq \frac{m}{n} \leq \sqrt{3D}$, where
\begin{equation}
\label{e_def}
e = \left\{ \begin{array}{ll}
0 & \mbox{if either $2 \mid D$, or $2 \mid (D+1), mn$} \\
1 & \mbox{otherwise.}
\end{array}
\right.
\end{equation}
Then
\begin{equation}
\label{z_IWR_mn}
\zeta^m_{\IWR(D)}(s) = \sum_{\substack{m,n \in \zed_{>0},\ \gcd(m,n)=1 \\ \sqrt{\frac{D}{3}} \leq \frac{m}{n} \leq \sqrt{3D}}} \left( \frac{2^e \gcd(m,D)}{m^2 + Dn^2} \right)^s,
\end{equation}
and so for each real $s > 1$,
\begin{eqnarray}
\label{z_IWR_mn_up}
\zeta^m_{\IWR(D)}(s) & \leq & (2D)^s  \sum_{\substack{m,n \in \zed \setminus \{0\} \\ \sqrt{\frac{D}{3}}  \leq \frac{m}{n} \leq \sqrt{3D}}} \frac{1}{ \left( m^2 + Dn^2 \right)^s} \nonumber \\
& \leq & \left( 2D \right)^s \sum_{m,n \in \zed \setminus \{0\}} \frac{1}{ \left( m^2 + Dn^2 \right)^s} =  \left( 2D \right)^s \zeta_{\que(\sqrt{-D})}(s).
\end{eqnarray}
Now, the Dedekind zeta-function of a number field converges absolutely for $\Re(s) > 1$ and has a simple pole at $s=1$.

On the other hand, for all real $s >1$,
\begin{eqnarray}
\label{z_IWR_mn_low}
\zeta^d_{\IWR(D)}(s) & \geq &  \sum_{\substack{m,n \in \zed \setminus \{0\},\ \gcd(m,n)=1 \\ \sqrt{\frac{D}{3}}  \leq \frac{m}{n} \leq \sqrt{3D}}} \frac{1}{\left( 2mn \right)^s}  \nonumber \\
& \geq & \frac{1}{\left( 2\sqrt{3D} \right)^s} \sum_{n=1}^{\infty} \frac{a_n}{n^{2s}},
\end{eqnarray}
where $a_n$ is the cardinality of the set
$$S_n = \left\{ m \in \zed_{>0} : n\sqrt{\frac{D}{3}} \leq m \leq n \sqrt{3D},\ \gcd(m,n)=1 \right\}.$$
We will now produce a lower bound on $a_n$ for every $n \geq 1$. For each $m \in S_n$, let $s_n(m) = m \md n$, then
$$a_n = |S_n| \geq \left| \left\{ s_n(m) : m \in S_n \right\} \right|.$$
Notice that 
$$\sqrt{3D} - \sqrt{D/3} = \sqrt{D} (\sqrt{3} - 1/\sqrt{3}) > 1$$
for each $D$, and hence
$$\left\{ s_n(m) : m \in S_n \right\} = \left\{ k \in \zed :  1 \leq k < n, \gcd(k,n) =1 \right\},$$
meaning that $a_n \geq \varphi(n)$, the Euler $\varphi$-function of $n$. Therefore
\begin{equation}
\label{zeta_phi}
\zeta^d_{\IWR(D)}(s) \geq \frac{1}{\left( 2\sqrt{3D} \right)^s} \sum_{n=1}^{\infty} \frac{\varphi(n)}{n^{2s}} = \frac{1}{\left( 2\sqrt{3D} \right)^s} \frac{\zeta(2s-1)}{\zeta(2s)}
\end{equation}
for all real $s > 1$ by Theorem 288 of \cite{hardy}. The right hand side of \eqref{zeta_phi} converges absolutely for  $\Re(s) > 1$ and has a simple pole at $s=1$. The inequality \eqref{zeta_bnd} now follows upon combining \eqref{zeta_md_ineq} with \eqref{z_IWR_mn_up} and~\eqref{zeta_phi}.

Since each Dirichlet series can be written in the form $\sum_{n=1} b_n n^{-s}$ for some coefficient sequence $\{b_n\}_{n=1}^{\infty}$, we will refer to $\sum_{n=1}^N b_n$ as its $N$-th partial sum of coefficients. Now Theorem~4.20 of \cite{sautoy} guarantees that the $N$-th partial sums of coefficients of Dirichlet series $\frac{1}{\left( 2\sqrt{3D} \right)^s} \frac{\zeta(2s-1)}{\zeta(2s)}$ and $\left( \frac{4D}{\sqrt{3}} \right)^s \zeta_{\que(\sqrt{-D})}(s)$ are equal to $O(N)$ as $N \to \infty$. Inequality \eqref{zeta_bnd} implies that the same must be true about $N$-th partial sums of coefficients of Dirichlet series $\zeta^m_{\IWR(D)}(s)$ and $\zeta^d_{\IWR(D)}(s)$, and that $\zeta^m_{\IWR(D)}(s)$ and $\zeta^d_{\IWR(D)}(s)$ are absolutely convergent for $\Re(s) > 1$ with limits in \eqref{lim_md} existing and nonzero for $s \in \real$. This finishes the proof of the theorem.
\endproof

\begin{rem} \label{ht_zeta} There is a connection between the zeta-function $\zeta^m_{\IWR(D)}(s)$ and the height zeta-function of the corresponding Pell-type rational conic. One can define a height function on points $\bx = (x_1, x_2,x_3) \in \zed^3$ as
$$H(\bx) = \frac{1}{\gcd(x_1,x_2,x_3)} \max_{1 \leq i \leq 3} |x_i|.$$
It is easy to see that $H$ is in fact projectively defined, and hence induces a function on a rational projective space. Let $D$ be a fixed positive squarefree integer, then the set of all integral points $(p,r,q)$ satisfying
\begin{equation}
\label{int_con}
p^2+Dr^2=q^2,\ \gcd(p,r,q)=1,\ q > 0
\end{equation}
is precisely the set of all distinct representatives of projective rational points on the Pell-type conic
$$X_D(\que) = \{ [x,y,z] \in \pee(\que^3) : x^2+Dy^2=z^2 \}.$$
For each point $[x,y,z] \in X_D(\que)$ there is a unique $(p,r,q)$ satisfying \eqref{int_con}, and 
$$H([x,y,z]) = H(p,r,q)=q.$$ 
Hence the height zeta-function of $X_D(\que)$ is
$$\sum_{[x,y,z] \in X_D(\que)} \frac{1}{H([x,y,z])^s} = \sum_{(p,r,q) \text{ as in \eqref{int_con}}} \frac{1}{q^s},$$
where $s \in \cee$.
\end{rem}
\bigskip

\section{IWR sublattices of IWR lattices}
\label{IWR_subl}

In this section we further investigate distribution properties of planar IWR lattices and prove Theorem~\ref{main}. Theorem~1.3 of \cite{fletcher_jones} guarantees that every IWR lattice of type $D$ contains IWR sublattices belonging to every similarity class of this type, and none others. Hence $\zeta^m_{\IWR(D)}(s)$ and $\zeta^d_{\IWR(D)}(s)$ are zeta-functions of minimal lattices over similarity classes of IWR sublattices of any IWR lattice of type $D$ in the plane. It will be convenient to define 
$$\Omega_D(p,q) = \sqrt{q}\ \Gamma_D(p,q) = \begin{pmatrix} q & p \\ 0 & r\sqrt{D} \end{pmatrix} \zed^2$$
for each $(p,r,q,D)$ satisfying \eqref{prqD}. Then for a fixed choice of $D,p_0,q_0$ the lattice $\Omega_D(p_0,q_0)$ contains IWR sublattices similar to each $\Omega_D(p,q)$. We will now  describe explicitly how these sublattices look like. We start with a simple example of such lattices.

\begin{lem} \label{two_IWR} Let $(p,r,q,D)$ and $(p_0,r_0,q_0,D)$ satisfy \eqref{prqD}. Let 
\begin{equation}
\label{kmn}
k=m^2+Dn^2
\end{equation}
for some $m,n \in \zed$, not both zero, and let 
\begin{equation}
\label{U}
U = \begin{pmatrix} \frac{m}{\sqrt{k}} & -\frac{n\sqrt{D}}{\sqrt{k}} \\ \frac{n\sqrt{D}}{\sqrt{k}} & \frac{m}{\sqrt{k}} \end{pmatrix}.
\end{equation}
Then $U$ is a real orthogonal matrix such that the lattice
$$\Lambda = \sqrt{k}\ r_0q_0 U \Omega_D(p,q)$$
is an IWR sublattice of $\Omega_D(p_0,q_0)$ similar to $\Omega_D(p,q)$ with
$$\left| \Omega_D(p_0,q_0) : \Lambda \right| = r_0q_0rqk.$$
\end{lem}

\proof As indicated in the proof of Theorem~1.3 of \cite{fletcher_jones},
$$\begin{pmatrix} q_0 & p_0 \\ 0 & r_0\sqrt{D} \end{pmatrix} \begin{pmatrix} r_0q & r_0p-rp_0 \\ 0 & rq_0 \end{pmatrix} = r_0q_0 \begin{pmatrix} q & p \\ 0 & r\sqrt{D} \end{pmatrix},$$
and so $r_0q_0 \Omega_D(p,q)$ is a sublattice of $\Omega_D(p_0,q_0)$ of index $rq$. Now notice that
\begin{eqnarray*}
\Lambda & = & \sqrt{k} r_0q_0 \begin{pmatrix} \frac{m}{\sqrt{k}} & -\frac{n\sqrt{D}}{\sqrt{k}} \\ \frac{n\sqrt{D}}{\sqrt{k}} & \frac{m}{\sqrt{k}} \end{pmatrix} \begin{pmatrix} q & p \\ 0 & r\sqrt{D} \end{pmatrix} \zed^2 \\
& = & \begin{pmatrix} q_0 & p_0 \\ 0 & r_0\sqrt{D} \end{pmatrix} \begin{pmatrix} mr_0q - np_0q & m(r_0p-p_0r)-n(p_0p+Dr_0r) \\ nq_0q & mrq_0+npq_0 \end{pmatrix} \zed^2,
\end{eqnarray*}
and hence is a sublattice of $\Omega_D(p_0,q_0)$ of index $r_0q_0rqk$.
\endproof

Lemma~\ref{two_IWR} demonstrates some examples of sublattices of $\Omega_D(p_0,q_0)$ similar to $\Omega_D(p,q)$. We will now describe all such sublattices.

\begin{thm} \label{two_IWR_1} A sublattice $\Lambda$ of $\Omega_D(p_0,q_0)$ is similar to $\Omega_D(p,q)$ as above if and only if
\begin{equation}
\label{subl}
\Lambda = \sqrt{Q_{p_0,q_0,p,q}(m,n)}\ U \Gamma_D(p,q),
\end{equation}
for some $m,n \in \zed$, not both zero, where $Q_{p_0,q_0,p,q}(m,n)$ is a positive definite binary quadratic form, given by \eqref{Q} below, and $U$ is a real orthogonal matrix as in \eqref{U_orth} with the angle $t$ satisfying \eqref{sin_cos}, where $x,y$ are as in \eqref{cong_sol_1} or  \eqref{cong_sol_2}. In this case,
\begin{equation}
\label{index}
\left| \Omega_D(p_0,q_0) : \Lambda \right| = \frac{rQ_{p_0,q_0,p,q}(m,n)}{r_0q_0}.
\end{equation}
\end{thm}

\proof  By Theorem~1.1 of \cite{fletcher_jones}, $\Lambda \sim \Omega_D(p,q)$ if and only if
\begin{equation}
\label{L1}
\Lambda = \sqrt{\frac{k}{q}}\ U  \begin{pmatrix} q & p \\ 0 & r\sqrt{D} \end{pmatrix} \zed^2
\end{equation}
for some positive integer $k$ and a real orthogonal matrix
\begin{equation}
\label{U_orth}
U = \begin{pmatrix} \cos t & -\sin t \\ \sin t & \cos t \end{pmatrix} \text{ or } \begin{pmatrix} \cos t & \sin t \\ \sin t & -\cos t \end{pmatrix}
\end{equation}
for some value of the angle $t$. On the other hand, $\Lambda \subset \Omega_D(p_0,q_0)$ if and only if
\begin{equation}
\label{L2}
\Lambda = \begin{pmatrix} q_0 & p_0 \\ 0 & r_0\sqrt{D} \end{pmatrix} C \zed^2,
\end{equation}
where $C$ is an integer matrix. Therefore $\Lambda$ as in \eqref{L1} is a sublattice of $\Omega_D(p_0,q_0)$ if and only if it is of the form \eqref{L2} with
$$C = \alpha\ \begin{pmatrix} q(r_0\sqrt{D}\cos t - p_0\sin t) & (r_0p-rp_0)\sqrt{D}\cos t - (pp_0+rr_0D)\sin t \\ qq_0\sin t & q_0p\sin t + q_0r\sqrt{D}\cos t \end{pmatrix}$$
or
$$C = \alpha\ \begin{pmatrix} q(r_0\sqrt{D}\cos t - p_0\sin t) & (r_0p+rp_0)\sqrt{D}\cos t - (pp_0-rr_0D)\sin t \\ qq_0\sin t & q_0p\sin t - q_0r\sqrt{D}\cos t \end{pmatrix}$$
where $\alpha = \frac{\sqrt{k}}{q_0r_0 \sqrt{qD}}$. These conditions imply that we must have
\begin{equation}
\label{sin_cos}
\cos t =  \frac{xp_0+yq_0}{\sqrt{qk}},\ \sin t = \frac{xr_0\sqrt{D}}{\sqrt{qk}}
\end{equation}
for some integers $x,y$ satisfying one of the following two systems of congruences:
\begin{equation}
\label{cong_1}
\left. \begin{array}{ll}
q_0rx + (p_0r - r_0p)y \equiv 0 (\md qr_0) \\
(p_0r + r_0p)x + q_0ry \equiv 0 (\md qr_0)
\end{array}
\right\},
\end{equation}
or
\begin{equation}
\label{cong_2}
\left. \begin{array}{ll}
q_0rx + (p_0r + r_0p)y \equiv 0 (\md qr_0) \\
(p_0r - r_0p)x + q_0ry \equiv 0 (\md qr_0)
\end{array}
\right\}.
\end{equation}

First assume \eqref{cong_1} is satisfied. Notice that
$$\det \begin{pmatrix} q_0r & p_0r-r_0p \\ p_0r+r_0p & q_0r \end{pmatrix} = (qr_0)^2 \equiv 0 (\md qr_0),$$
which means that a pair $(x,y)$ solves the system \eqref{cong_1} if and only if it solves one of these two congruences. Hence it is enough to solve the first congruence of  \eqref{cong_1}. Define $d_1 = \gcd(q_0r,qr_0)$ and $d_2 = \gcd(d_1,p_0r - r_0p)$, and let $a,b \in \zed$ be such that
$$aq_0r+bqr_0=d_1.$$
It now easily follows that the set of all possible solutions to \eqref{cong_1} is
\begin{equation}
\label{cong_sol_1}
(x,y) = \left\{ \left( \frac{a(r_0p-p_0r)n}{d_2} + \frac{qr_0m}{d_1},\ \frac{d_1n}{d_2} \right) : n,m \in \zed \right\}.
\end{equation}
Combining \eqref{sin_cos} with \eqref{cong_sol_1}, we see that
\begin{equation}
\label{qk1}
qk = \left( \frac{\left( ap_0(r_0p-p_0r) + d_1q_0 \right)n}{d_2} + \frac{qp_0r_0m}{d_1} \right)^2 + Dr_0^2 \left( \frac{a(r_0p-p_0r)n}{d_2} + \frac{qr_0m}{d_1} \right)^2.
\end{equation}
Then the right hand side of \eqref{qk1} is a positive definite integral binary quadratic form in the variables $m,n$:
\begin{eqnarray}
\label{Q1}
Q^1_{p_0,q_0,p,q}(m,n) & = & \left\{ \frac{q_0^2 q^2 r_0^2}{d_1^2} \right\} \times m^2 \nonumber \\
& + & \left\{ \frac{ a^2 (r_0p-p_0r)^2 q_0^2 + 2ad_1p_0q_0(r_0p-p_0r) + d_1^2q_0^2}{d_2^2} \right\} \times n^2 \nonumber \\
& + & \left\{ \frac{2a(r_0p-p_0r)qq_0^2r_0 + 2d_1qq_0r_0p_0}{d_1d_2} \right\} \times mn.
\end{eqnarray}
One can observe that all three coefficients of $Q^1_{p_0,q_0,p,q}(m,n)$ are divisible by $q$. Then define
\begin{equation}
\label{Q}
Q_{p_0,q_0,p,q}(m,n) = \frac{1}{q} Q^1_{p_0,q_0,p,q}(m,n),
\end{equation}
which is again a positive definite integral binary quadratic form.

Now notice that the system of congruences in \eqref{cong_2} is the same as the one in \eqref{cong_1} with the order of equations reversed and the variables $x$ and $y$ reversed. Hence the solution set for \eqref{cong_2} is
\begin{equation}
\label{cong_sol_2}
(x,y) = \left\{ \left( \frac{d_1n}{d_2},\ \frac{a(r_0p-p_0r)n}{d_2} + \frac{qr_0m}{d_1} \right) : n,m \in \zed \right\}.
\end{equation}
Combining \eqref{sin_cos} with \eqref{cong_sol_2}, we see that if \eqref{cong_2} is satisfied, then
\begin{equation}
\label{qk2}
qk = \left( \frac{\left( aq_0(r_0p-p_0r)+d_1p_0 \right)n}{d_2} + \frac{qq_0r_0m}{d_1} \right)^2 + \frac{Dr_0^2d_1^2n^2}{d_2^2}.
\end{equation}
Then the right hand side of \eqref{qk2} is precisely $Q^1_{p_0,q_0,p,q}(n,m)$. 

In either case, we have
\begin{equation}
\label{k_value}
k = \frac{1}{q} Q^1_{p_0,q_0,p,q}(m,n) = Q_{p_0,q_0,p,q}(m,n)
\end{equation}
for some $m,n \in \zed$, not both zero. Then \eqref{subl} follows upon combining \eqref{L1} with \eqref{k_value}. Now we notice that
$$\left| \Omega_D(p_0,q_0) : \Lambda \right| = \frac{\det \Lambda}{\det \Omega_D(p_0,q_0)},$$
and so \eqref{index} follows from \eqref{subl}. This completes the proof of the lemma.
\endproof

Now define $S_D(p_0,q_0)$ to be the set of all IWR sublattices of $\Omega_D(p_0,q_0)$, and $S_D(p_0,q_0,p,q)$ to be the set of all IWR sublattices of $\Omega_D(p_0,q_0)$ which are similar to $\Omega_D(p,q)$. Then
$$S_D(p_0,q_0) = \bigsqcup S_D(p_0,q_0,p,q).$$
Define
\begin{equation}
\label{Z_D_1}
Z_{D,p_0,q_0,p,q}(s) = \sum_{\Lambda \in S_D(p_0,q_0,p,q)} \frac{1}{\left| \Omega_D(p_0,q_0) : \Lambda \right|^s}
\end{equation}
and
\begin{equation}
\label{Z_D_2}
Z_{D,p_0,q_0}(s) = \sum_{\Lambda \in S_D(p_0,q_0)} \frac{1}{\left| \Omega_D(p_0,q_0) : \Lambda \right|^s} = \sum_{(p,q) \text{ as in \eqref{prqD}}} Z_{D,p_0,q_0,p,q}(s)
\end{equation}
for $s \in \cee$.

\begin{lem} \label{zeta_two} For every squarefree positive integer $D$ and integer triples $(p_0,r_0,q_0)$ and $(p,r,q)$ satisfying \eqref{prqD}, the Dirichlet series $Z_{D,p_0,q_0,p,q}(s)$ is absolutely convergent for $\Re(s) > 1$. Moreover, it has analytic continuation to all of $\cee$ except for a simple pole at $s=1$.
\end{lem}

\proof By Theorem~\ref{two_IWR_1},
\begin{equation}
\label{eps1}
Z_{D,p_0,q_0,p,q}(s) = \left( \frac{r_0q_0}{r} \right)^s \sum_{(m,n) \in \zed^2 \setminus \{ \bo \}} \frac{1}{Q_{p_0,q_0,p,q}(m,n)^s},
\end{equation}
where the sum on the right hand side of \eqref{eps1} is the Epstein zeta-function of the positive definite integral binary quadratic form $Q_{p_0,q_0,p,q}(m,n)$; it is known to converge absolutely for $\Re(s) > 1$ and has analytic continuation to all of $\cee$ except for a simple pole at $s=1$ (this is a classical result, which can be found for instance in Chapter~5, \S 5 of \cite{koecher}; in fact, the authors of \cite{koecher} indicate that the existence of a simple pole at $s=1$ goes as far back as the work of Kronecker, 1889). The lemma follows.
\endproof

\begin{lem} \label{zeta_two_1} For every squarefree positive integer $D$ and integer triple $(p_0,r_0,q_0)$ satisfying \eqref{prqD}, the Dirichlet series $Z_{D,p_0,q_0}(s)$ is absolutely convergent for $\Re(s) > 1$ and for $s \in \real$ the limit
\begin{equation}
\label{lim_ZD}
\lim_{s \to 1^+} (s-1)^2 Z_{D,p_0,q_0}(s)
\end{equation}
exists and is nonzero. Moreover, if we write $Z_{D,p_0,q_0}(s) = \sum_{n=1}^{\infty} b_n n^{-s}$, then the $N$-th partial sum of coefficients of $Z_{D,p_0,q_0}(s)$ is
$$\sum_{n=1}^N b_n = O(N \log N)$$
as $N \to \infty$.
\end{lem}

\proof Combining \eqref{Z_D_2}, \eqref{eps1}, and \eqref{rq} we obtain for every real $s>0$
$$Z_{D,p_0,q_0}(s) = (r_0q_0)^s \sum_{(p,r,q) \text{ as in \eqref{prqD}}} \left( \frac{1}{r^s} \sum_{(m,n) \in \zed^2 \setminus \{ \bo \}} \frac{1}{Q_{p_0,q_0,p,q}(m,n)^s} \right).$$
Combining this observation with Theorem~\ref{IWR_zeta} implies that
\begin{eqnarray}
\label{eps2}
&\ & \left( \frac{r_0q_0}{2\sqrt{3}} \right)^s \frac{\zeta(2s-1)}{\zeta(2s)} \inf_{(p,r,q) \text{ as in \eqref{prqD}}} \sum_{(m,n) \in \zed^2 \setminus \{ \bo \}} \frac{1}{Q_{p_0,q_0,p,q}(m,n)^s} \\
& \leq & \left( r_0q_0 \sqrt{D} \right)^s \zeta^d_{\IWR(D)}(s) \inf_{(p,r,q) \text{ as in \eqref{prqD}}} \sum_{(m,n) \in \zed^2 \setminus \{ \bo \}} \frac{1}{Q_{p_0,q_0,p,q}(m,n)^s} \nonumber \\
& \leq & Z_{D,p_0,q_0}(s) \nonumber \\
&\leq & \left( r_0q_0 \sqrt{D} \right)^s \zeta^d_{\IWR(D)}(s) \sup_{(p,r,q) \text{ as in \eqref{prqD}}} \sum_{(m,n) \in \zed^2 \setminus \{ \bo \}} \frac{1}{Q_{p_0,q_0,p,q}(m,n)^s} \nonumber \\
&\leq & \left( \frac{4r_0q_0D^{\frac{3}{2}}}{\sqrt{3}} \right)^s \zeta_{\que(\sqrt{-D})}(s) \sup_{(p,r,q) \text{ as in \eqref{prqD}}} \sum_{(m,n) \in \zed^2 \setminus \{ \bo \}} \frac{1}{Q_{p_0,q_0,p,q}(m,n)^s}. \nonumber
\end{eqnarray}
Theorem~\ref{IWR_zeta} and Lemma~\ref{zeta_two} now imply that for each $(p,r,q)$ as in \eqref{prqD} the Dirichlet series
\begin{equation}
\label{Dir1}
\left( \frac{r_0q_0}{2\sqrt{3}} \right)^s \frac{\zeta(2s-1)}{\zeta(2s)} \sum_{(m,n) \in \zed^2 \setminus \{ \bo \}} \frac{1}{Q_{p_0,q_0,p,q}(m,n)^s}
\end{equation}
and
\begin{equation}
\label{Dir2}
\left( \frac{4r_0q_0D^{\frac{3}{2}}}{\sqrt{3}} \right)^s \zeta_{\que(\sqrt{-D})}(s) \sum_{(m,n) \in \zed^2 \setminus \{ \bo \}} \frac{1}{Q_{p_0,q_0,p,q}(m,n)^s}
\end{equation}
are absolutely convergent for $\Re(s) > 1$ and have analytic continuation to the half-plane $\Re(s) > 0$ except for a pole of order 2 at $s=1$. Then Theorem~4.20 of \cite{sautoy} implies that the $N$-th partial sums of coefficients of all the Dirichlet series as in \eqref{Dir1} and \eqref{Dir2} must be equal to $O(N \log N)$. Then \eqref{eps2} implies that the $N$-th partial sum of coefficients of $Z_{D,p_0,q_0}(s)$ is also $O(N \log N)$, and $Z_{D,p_0,q_0}(s)$ is absolutely convergent for $\Re(s) > 1$, where the limit of~\eqref{lim_ZD} exist and is nonzero for $s \in \real$.
\endproof

\proof[Proof of Theorem~\ref{main}] The theorem now follows upon combining Lemmas~\ref{reduce},~\ref{sim_red}, and~\ref{zeta_two_1}.
\endproof
\bigskip

{\bf Acknowledgment.} I would like to thank Stefan K\"uhnlein for providing me with a copy of his manuscript \cite{kuehnlein}, as well as for discussing the problem with me, and for his useful remarks about this paper. I would also like to thank Michael D. O'Neill and Bogdan Petrenko for their helpful comments on this paper. Finally, I would like to thank the referee for the suggestions and corrections which improved the quality of this paper.
\bigskip

\bibliographystyle{plain}  
\bibliography{iwr_zeta}    
\end{document}